\documentclass[11pt,bezier]{article}
\usepackage{amsmath,graphicx,amssymb,amsfonts}
\renewcommand{\baselinestretch}{1.3}
\newtheorem{prethm}{{\bf Theorem}}

\newenvironment{thm}{\begin{prethm}{\hspace{-0.5
               em}{\bf.}}}{\end{prethm}}
\newtheorem{prepro}{{\bf Theorem}}

\newtheorem{precor}{{\bf Corollary}}

\newenvironment{cor}{\begin{precor}{\hspace{-0.5
               em}{\bf.}}}{\end{precor}}
\newtheorem{preconj}{{\bf Conjecture}}

\newtheorem{preremark}{{\bf Remark}}

\newenvironment{remark}{\begin{preremark}\em{\hspace{-0.5
               em}{\bf.}}}{\end{preremark}}
\newtheorem{prelem}{{\bf Lemma}}

\newtheorem{eg}{Example}
\newenvironment{example}{\begin{eg}\rm}{\end{eg}}
\newtheorem{preproof}{{\bf Proof.}}

\newenvironment{proof}[1]{\begin{preproof}{\rm
               #1}\hfill{$\Box$}}{\end{preproof}}

\title{\Large {\bf Inclusion Matrices and Chains}}

\author{
E. Ghorbani$^{\textrm{b,c}}$ \and  G.B.
Khosrovshahi$^{\textrm{c,a}}$ \and Ch. Maysoori$^{\textrm{c}}$ \and
M. Mohammad-Noori$^{\textrm{a},\textrm{c}}$ \\
{\small$^{\textrm{a}}$\it Center of Excellence, School of
Mathematics, Statistics, and Computer Science,}\\
 {\small\it University of Tehran, Tehran, Iran}\\
{\small$^{\textrm{b}}$\it Department of Mathematical Sciences,
 Sharif University of Technology, Tehran, Iran}\\
 {\small$^{\textrm{c}}$\it Institute for Studies in
Theoretical Physics and Mathematics (IPM),}\\
{\small\it P.O. Box 19395-5746, Tehran, Iran}\\
{\footnotesize{$\mathsf{e\_ghorbani@math.sharif.edu}$}\quad\quad
{$\mathsf{rezagbk@ipm.ir}$}\quad\quad
{$\mathsf{maysoori@ipm.ir}$}\quad\quad
{$\mathsf{morteza@ipm.ir}$}}
 }
\date{}
\begin{document}
\maketitle
\begin{quote}
{\small \hfill{\rule{13.3cm}{.1mm}\hskip2cm}
\textbf{Abstract}\vspace{1mm} {\renewcommand{\baselinestretch}{1}
\parskip = 0 mm

Given integers $t$, $k$, and $v$ such that $0\leq t\leq k\leq v$,
let $W_{tk}(v)$ be the inclusion matrix of $t$-subsets vs.
$k$-subsets of a $v$-set.
We modify slightly the concept of standard tableau to study the
notion of rank of a finite set of positive integers which was
introduced by Frankl. Utilizing this, a decomposition of the poset
$2^{[v]}$ into symmetric skipless chains is given. Based on this
decomposition, we construct an inclusion matrix, denoted by
$W_{\overline{t}k}(v)$, which is row-equivalent to $W_{tk}(v)$. Its
Smith normal form is determined. As applications, Wilson's
 diagonal form of $W_{tk}(v)$ is obtained as well as a new proof of
 the well known theorem on the necessary and sufficient conditions
for existence of integral solutions of the system
$W_{tk}\bf{x}=\bf{b}$ due to Wilson.
 Finally we present another inclusion matrix
  with similar properties to those of $W_{\overline{t}k}(v)$
 which is in some way equivalent to  $W_{tk}(v)$. }}

 \vspace{0.1mm}

\noindent{\small {\it 2000 Mathematics Subject Classification}:
05B05; 05B20; 15A21; 05D05

\vspace{-3mm}\hfill{\rule{13.3cm}{.1mm}\hskip2cm}}
\end{quote}
\newpage
\vspace{9mm} \noindent{\bf\Large 1. Introduction}\vspace{3mm}

\noindent To start we fix some notations.
For a finite set $F$ and a nonnegative integer $i$ let ${F \choose
i}$ be the set of all $i$-subsets of $F$, and  $2^F$ be the set of
all subsets of $F$. For a positive integer $v$, we denote the set
$\{1,\ldots,v\}$ by $[v]$ and define ${v\choose -1}$ to be 0.

Let $t, k,$ and $v$ be integers satisfying $0\leq t\leq k\leq v$.
The inclusion matrix $W_{tk}(v)$ (denoted briefly as $W_{tk}$) is a
$0,1$-matrix whose rows and columns are indexed  by $t$-subsets and
$k$-subsets of the set $[v]$, respectively, and $W_{tk}(v)(T,K)=1$
if and only if $T\subseteq K$. These inclusion matrices have
applications in many combinatorial problems, particularly in design
theory \cite{Beth} and their properties are studied extensively
\cite{fra, gju, kmHgh, wil2, wilson}.

Our motivation in this work relates to the simple but fundamental
equation
\begin{equation}\label{w1} W_{tk}{\bf x}=\lambda {\bf 1},\end{equation}
 where $\lambda$ is a positive integer and ${\bf 1}$
 is the all 1's vector. Any
 integral solution of (\ref{w1}) is called a
$t$-$(v,k,\lambda)$ signed design; and any non-negative integral
solution of (\ref{w1}) is a $t$-design. Therefore of importance is having an equivalent system to the above
such that its coefficient matrix is of a `simpler' nature.

 Our aim in this paper is to construct such matrices using a
 notion of  `rank' of a finite set. Our
 criterion to being simpler is the Smith form. The Smith form of
 the new matrices is roughly identity while that of $W_{tk}$ is more
 complicated. See Theorems \ref{snf}, \ref{last} and Corollary \ref{wil}.

The paper is organized as follows: In Section 2 we use a slightly
modified version of standard tableau to study the notion of the
`rank' of a finite set of positive integers, introduced by Frankl
\cite{fra}. This notion is useful in Section 3 to describe a
decomposition of $2^{[v]}$ into symmetric chains, in which the
elements of each chain have the same rank. This decomposition
coincides with that of De Bruijn described recursively in \cite{deb}
and explicitly in \cite{and, gk}. In Section 4 we introduce a new
inclusion matrix $W_{\overline{t}k}(v)$ based on the above
decomposition as follows: Replace each $t$-subset $T$ with the
minimal element $\overline{T}$ in the same chain and construct the
related inclusion matrix $W_{\overline{t}k}(v)$. The new matrix
$W_{\overline{t}k}(v)$ is row-wise equivalent to $W_{tk}$  and has
$(I\mid O)$ as the Smith normal form. This is implicit in the work
of Bier \cite{Bier}. Using these results we obtain Wilson's
 diagonal form of $W_{tk}(v)$ as well as a new proof of
 the well known theorem on the existence of signed (integral)
 $t$-designs due to Wilson and Graver and Jurkat.
 In Section 5 a similar approach is applied to find another
inclusion matrix $W_{t\underline{k}}$ which is somehow equivalent to
$W_{tk}(v)$ and have the same Smith form as $W_{\overline{t}k}(v)$.

\vspace{5mm} \noindent{\bf\Large 2. Rank of a finite set}
\vspace{3mm}

The rank of a given set $F$ of positive integers, denoted by $r(F)$,
is defined as follows \cite{fra}: To each $F\subseteq [v]$ we
associate a walk $w(F)$ from the origin to $(v-|F|,|F|)$ by steps of
length one such that the $i$-th step is to the right (up) if $i\not
\in F$ ($i \in F$); the rank of $F$, denoted as $r(F)$, is then
$|F|-b$, where $b$ is the largest integer such that the line $y=x+b$
touches $w(F)$ from the above. This definition does not depend on
$v$ in the sense that replacing $v$ by any integer $\geq \max(F)$,
the rank does not change.

A {\it standard tableau} is a $2\times m$ array, filled with
integers $1,\ldots,2m$, such that the entries across any row or down
any column are ordered increasingly. We modify this object to
provide a useful tool which reveals the  properties of rank:
Consider a $2 \times |F|$ tableau $\mathcal{T}(F)$ and arrange the
elements of $F$ on the first row in increasing order. Fill the
second row from left to right as follows: Below each element $a$ put
the largest positive integer smaller than $a$ not appearing on the
left neither in the second  row nor in the first; If such an integer
does not exist, put the symbol $\mathfrak{j}$. After filling the
second row, let $r(F)=|F|-\#\mathfrak{j}$'s. In other words
if$a_1<a_2<\cdots<a_{|F|}$ are the elements of $F$, then

$$\mathcal{T}(F)=\left(\begin{array}{rrrrr}
a_1 & a_2 & \ldots & a_{|F|}\\
b_1 & b_2& \ldots & b_{|F|}\end{array}\right).$$

\noindent In which, for $i=1,\ldots,|F|$, we have
\begin{equation}b_i=\max \left( [a_i]\setminus
(F\cup\{b_1,\ldots,b_{i-1}\}) \right),\end{equation}

\noindent where $\max \emptyset$ is defined to be $\mathfrak{j}$.
The set of non-$\mathfrak{j}$ elements in the second row is denoted
by fill$(F)$. The set of corresponding elements in the first row is
denoted by ${\rm fill}^*(F)$:
\begin{align*}
   {\rm fill} (F)&=\{b_i: i=1,2,\ldots, |F|\}\setminus \{\mathfrak{j}\},\\
{\rm fill}^* (F)&=\{a_i: b_i\in {\rm fill}(F)\}=\{a_i: b_i\neq
\mathfrak{j}\}.
\end{align*}
Hence
\begin{equation}\label{rfill}r(F)=|{\rm fill}(F)|=|{\rm fill}^*(F)|.\end{equation}

\noindent If $r(F)=|F|$, then $F$ is said to be {\it full-rank}.

\begin{example}\label{S2378}
Let $F=\{2,3,7,8\}$. Then
$$\nonumber\mathcal{T}(F)=\left(\begin{array}{rrrr}
2 & 3 &7 &8\\
1 & \mathfrak{j}& 6 & 5\end{array}\right).$$ Thus $r(F)=3$ and $F$
is not full-rank. One can consider the walk $w(F)$ corresponding $F$
to calculate the rank of $F$: $w(F)=RUURRUU$ and the corresponding
tangent line is $y=x+1$, so $r(F)=4-1=3$.
\end{example}

\vspace{5mm} \noindent{\bf\Large 3. A decomposition of the poset
$2^{[v]}$}\vspace{3mm}

In this section, first we introduce an algorithm to construct a
chain of sets of positive integers in which all the sets have the
same rank. Then we develop this method to decompose $2^{[v]}$ into
symmetric skipless chains for a given value of $v$. The algorithm is
extendable in the sense that applying it to $2^{[v+1]}$ and
eliminating $v+1$ yields the same result as applying it to
$2^{[v]}$.

To construct a chain which includes a given set $F$ of positive
integers, two methods are required: A method to find the successor
$F^+$ of $F$ and another one is to find the predecessor $F^-$ of
$F$. We  simply define  $F^+:=F\cup \{a\}$, where $a$ is the minimal
positive integer which has not appeared in $\mathcal{T}(F)$. On the
other hand, $F^-=F\setminus \{b\}$ where $b$ is the element above
the rightmost $\mathfrak{j}$ in $\mathcal{T}(F)$; If there is no
$\mathfrak{j}$ in the second row of $\mathcal{T}(F)$ then $F$ has
full rank and no predecessor. It is easy to prove that $(F^+)^-=F$
and if $F$ does not have full rank, then $(F^-)^+=F$.

\begin{example}\label{PM}
Consider the set $F$ and the corresponding tableau
$\mathcal{T}(F)$ as in Example \ref{S2378}: The number $4$ is the
least positive integer which has not appeared in ${\mathcal
T}(F)$, hence $F^{+}=F\cup\{4\}$. Now $\mathcal{T}(F^{+})$ is as
follows:

$$\nonumber \mathcal{T}(F^{+})=\left(\begin{array}{rrrrr}
2 & 3 &4 &7 &8\\
1 & \mathfrak{j} & \mathfrak{j}  &6 &5\end{array}\right).$$

\noindent It is immediately seen that ${\rm fill}(F^+)={\rm
fill}(F),\, {\rm fill}^*(F^+)={\rm fill}^*(F)$ so by (\ref{rfill}),
$r(F)=r(F^+)$. In fact the element below $4$ is $\mathfrak{j}$ and
adding $4$ to $F$, has no effect on how the second row is filled.
Thus the rank does not change. Given the set $E=F^+$, to recover $F$
consider the tableau $\mathcal{T}(E)$ and eliminate the element
above the rightmost $\mathfrak{j}$ from $E$ to obtain $E^-$. It is
seen that $F=E^-$.
\end{example}

Now we use the above algorithm to find a decomposition of $2^{[v]}$
into symmetric skipless chains for a fixed positive integer $v$.
Given a set $F\subseteq [v]$, we would like to construct the chain
which contains this subset. Let $q=|F|, p=r(F)$ and $A_q=F$. For $p
\leq i<q$, define $A_i=A_{i+1}^-$ and for $q\leq i<v-p$, let
$A_{i+1}=A_i^+$. The set $A_p$ has full rank so $A_p^-$ is not
defined; And, $A_p$ is the only element of the chain with this
property. On the other hand, $\mathcal{T}(A_{v-p})$ contains totally
$(v-p)+p=v$ integers, which means that $A_{v-p}^+$ is no more a
subset of $[v]$; Moreover, $A_{v-p}$ is the only element of the
chain with this property. Hence we have the following chain
$$A_p\to\cdots\to A_q \to\cdots\to A_{v-p}$$
\noindent in which $A_{i+1}=A_i^+$ for $p\leq i<v-p$. Now choose an
element $E$ which has not appeared in this chain and similarly
construct the corresponding chain. Continue in this way until no
more subset of $[v]$ is remained. Obviously by this construction
every element of $2^{[v]}$ appears in a unique chain.
 Each chain is determined by its first element which is a full-rank subset of $[v]$
  and denoted as $\overline{F}$ where $F$ is any set of the chain.

\begin{example}\label{chains}
For $v=6$, the chains are as follow: The only chain of rank $0$ is
{\footnotesize $$\emptyset \to 1 \to 12 \to 123 \to 1234 \to 12345
\to 123456. $$}
 Chains with rank $1$ are {\footnotesize
$$ 2 \to 23 \to 234 \to 2345 \to 23456,$$
$$ 3 \to 13 \to 134 \to 1345 \to 13456,$$
$$ 4 \to 14 \to 124 \to 1245 \to 12456,$$
$$ 5 \to 15 \to 125 \to 1235 \to 12356,$$
$$ 6 \to 16 \to 126 \to 1236 \to 12346.$$}
Chains with rank $2$ are {\footnotesize
$$ 24 \to 245 \to 2456,\ \ \ \ \ 25 \to 235 \to 2356, $$
$$26 \to 236 \to 2346,  \ \ \ \ \ 34 \to 345 \to 3456, $$
$$35 \to 135 \to 1356, \ \ \ \ \ 36 \to 136 \to 1346, $$
$$45 \to 145 \to 1456, \ \ \ \ \ 46 \to 146 \to 1246, $$
$$56 \to 156 \to 1256.$$}
 Finally, chains with rank $3$ are those with only one element
 as the following
{\footnotesize $$246, \ \ \ \ 256,\ \ \ \ 346,\ \ \ \ 356, \ \ \ \
456.$$}
\end{example}

\begin{remark}
Considering Example \ref{chains}, given a set $F$ and a positive
integer $m$, it is useful to have a straightforward method to obtain
elements in the same chain with distance $m$ from $F$. Let $F^{+m}$
(resp. $F^{-m}$) be the set in the same chain as $F$ which has $m$
more (resp. $m$ less) elements. Let $a_1<\cdots<a_{|F|}$ be the
elements of $F$, and
$$B=\cup_{i=1}^{|F|+1}\{x\in\mathbb{Z} : a_{i-1}+1\le x\le
a_i-2\},$$ where $a_0=0$ and $a_{|F|+1}=+\infty $. By induction on
$m$ it is seen that to construct $F^{+m}$ from $F$, it is enough to
add $m$ least elements of $B$ to $F$. To construct $F^{-m}$ from
$F$, one should delete the elements corresponding to $m$ rightmost
$\mathfrak{j}$'s in $\mathcal{T}(F)$ from $F$.
\end{remark}

\begin{remark}\label{i-sub} The number of full-rank subsets of $[v]$ with rank $\leq r$, for
$0\leq r\leq v/2$, is ${v \choose r}$. To see this
 it is enough to correspond to each subset $F$ with
$r(F)\leq r$, the $r$-subset in the same chain. Then it follows that
the number of full-rank subsets with rank $r$
 (or equivalently the number of
chains with rank $r$) is ${v \choose r}-{v \choose {r-1}}$.
\end{remark}

 \vspace{5mm} \noindent{\bf\Large 4. The inclusion matrix
$W_{\overline{t}k}(v)$}\vspace{3mm}

We begin this section by recalling some definitions from matrix
theory. A {\it unimodular matrix} is a square integral matrix with
determinant $\pm1$, or equivalently whose inverse is integral. It is
well known that for a given integral matrix $A$ of rank $r$, there
exist unimodular matrices $U$ and $V$ such that $UAV=(D\mid O)$, in
which $D$ is a diagonal matrix with the positive integers $d_1,
\ldots ,d_r$ on its diagonal such that $d_1|\cdots |d_r$. The matrix
$(D\mid O)$ is called {\it Smith normal form} of $A$. Moreover,
\begin{equation}\label{di}d_i=f_i/f_{i-1},\end{equation} where $f_0=1$ and
$f_k$ is the greatest
common divisor of all minors of $A$ of order $k$, $1\leq k\leq r$.
For more on Smith form see \cite[pp. 26--33]{new}.

In this section the rank chains which was introduced in the previous
section is used to construct a new inclusion matrix from $W_{tk}$.
Any row index $T$ of $W_{tk}$ is replaced by $\overline{T}$,  the
indices of columns are kept, and the new 0,1-entries, like those of
$W_{tk}$,
  are determined by the inclusion relation. We denote this
matrix by $W_{\overline{t}k}$. We shall prove that the row space of
$W_{\overline{t}k}(v)$ is the same as that of $W_{tk}(v)$. It is
also shown that the Smith form of $W_{\overline{t}k}$ is $(I\mid
O)$, where $I$ is the identity matrix.

Let $R_{it}$ be the inclusion matrix whose rows are indexed by the
all $i$-subsets of $[v]$ of rank $i$, and the columns by $[v]\choose
t$. Note that by Remark {\ref{i-sub}}, $R_{it}$ has exactly
${v\choose i}-{v\choose i-1}$ rows. Then
\begin{equation}\label{wr1}W_{\overline{t}k}={\small\begin{array}{|c|}
                           \hline  R_{0k} \\\hline
                           R_{1k} \\\hline
                           \vdots \\\hline
                           R_{tk} \\\hline
                            \end{array}}\,.\end{equation}
 We observe that \begin{equation} \label{rw} R_{it}W_{tk}= {k-i\choose t-i}R_{ik}
,\quad\quad i\leq t\leq k \leq v.\end{equation} This holds because
for an $i$-subset $S$ and a $k$-subset $K$ of $[v]$,
$$(R_{it}W_{tk})(S,K)=\sum_{T\in{[v]\choose t}}R_{it}(S,T)W_{tk}(T,K).$$
 The right-hand side is the number of $t$-subsets $T$ such that $S\subseteq T\subseteq
 K$, and this number is ${k-i\choose t-i}$ if $S\subseteq K$, and 0
 otherwise.

It follows from (\ref{rw}) that
\begin{equation}\label{wr2}W_{\overline{i}t}W_{tk}=
                            {\small\begin{array}{|c|}
                           \hline {k\choose t}R_{0k} \\\hline
                           {k-1\choose t-1}R_{1k} \\\hline
                           \vdots \\\hline
                           {k-i\choose t-i}R_{ik} \\\hline
                           \end{array}}\,.
                        \end{equation}
Define the  ${v\choose t}\times{v\choose t}$ diagonal matrix
$D_{\overline{t}k}(v)$ to be  {\small$${\rm diag}\left({k-i \choose
t-i}^{{v\choose i}-{v\choose i-1}}, \quad i=0, 1, \ldots,
t\right),$$}  where the exponents indicates the multiplicity.  Then
by (\ref{wr2}),
\begin{equation}\label{wd}W_{\overline{t}t}W_{tk}=D_{\overline{t}k}W_{\overline{t}k}.
                        \end{equation}

\begin{thm}\label{snf} The Smith normal form of $W_{\overline{t}k}$, $t\le k \le v-t$, is
$(I\mid O)$, where $I$ is the identity matrix of order $v\choose t$.
\end{thm}
\begin{proof}{By (\ref{di}), it is enough to show that $W_{\overline{t}k}$ has a
unimodular submatrix of order $v\choose t$. By induction on $v+t$ we
show that the submatrix of $W_{\overline{t}k}(v)$ consists of the
columns indexed by a subset of rank $\leq t$, denoted by
$A_{t,k}(v)$, is the desired one. The assertion is clear for
$v\leq2$ with $t\leq1$. Let $v\geq3$, and $t\geq1$. First let $k=t$.
In an appropriate ordering we have
  \begin{equation*}\label{wtt}W_{\overline{t}t}={\small \begin{array}{|c|c|}
                           \hline A_{t-1,t} & * \\\hline
                           O & I \\\hline
                            \end{array}}\,,\end{equation*}
where $I$ is the identity matrix of order ${v\choose t}-{v\choose
t-1}$, whose rows and columns are indexed by subsets of rank $t$. By
induction, $A_{t-1,t}$ is unimodular and so is $W_{\overline{t}t}$.

Let $t < k < v-k$. The rows and columns of $W_{\overline{t}k}$
partitioned according to whether or not they contain the element $v$
as well as the following fact {\small
\begin{equation}\label{T}
   \{T\subseteq [v] : r(T)\leq t\}=\{T\cup\{v\}
: T\subseteq [v-1],\, r(T)\leq t-1\}\cup\{T\subseteq [v-1] :
r(T)\leq t\}.
\end{equation}
}
 This gives us the
following decomposition of $W_{\overline{t}k}$
\begin{equation}\label{part} W_{\overline{t}k}(v)={\small\begin{array}{|c|c|}
                           \hline W_{t\overline{-}1,k-1}(v-1) & O \\\hline
                           W_{\overline{t},k-1}(v-1) & W_{\overline{t}k}(v-1) \\\hline
                            \end{array}}\,.\end{equation}
Let $t<k<v-t$. By (\ref{T}), we observe that
$$ A_{t,k}(v)={\small\begin{array}{|c|c|}
                           \hline A_{t-1,k-1}(v-1) & O \\\hline
                           * & A_{t,k}(v-1) \\\hline\end{array}}\,.$$
 The induction hypothesis implies that the blocks on the diagonal are
unimodular, thus
$A_{t,k}(v)$ is unimodular.  Now let $k=v-t$.  Note that in this case the number of rows of the matrix $W_{\overline{t}k}(v-1)$
 is more than the number of its columns, so the Equation (\ref{part}) cannot be used directly.
 But using the partitions
$$W_{\overline{t}k}(v-1)={\small\begin{array}{|c|}
                           \hline W_{t\overline{-}1,k}(v-1) \\\hline
                            R_{tk}(v-1)  \\\hline
                            \end{array}\,,
                            \quad\quad
                           W_{\overline{t},k-1}(v-1)}={\small \begin{array}{|c|}
                           \hline W_{t\overline{-}1, k-1}(v-1) \\\hline
                            R_{t, k-1}(v-1) \\\hline
                            \end{array}}\,, $$
Equation (\ref{part}) becomes $$ W_{\overline{t}k}(v)={\small
\begin{array}{|c|c|}
                           \hline W_{t\overline{-}1,k-1}(v-1) & O \\\hline
                           W_{t\overline{-}1,k-1}(v-1) & W_{t\overline{-}1,k}(v-1) \\\hline
                           R_{t, k-1}(v-1)& R_{tk}(v-1)\\\hline \end{array}}\,.$$
Therefore it is easily seen that \begin{align*}\det
W_{\overline{t}k}(v)&= \pm\det\, {\small\begin{array}{|c|c|}\hline
W_{t\overline{-}1,k-1}(v-1) & O \\\hline
                            R_{t, k-1}(v-1)& R_{tk}(v-1)\\\hline
                            O & W_{t\overline{-}1,k}(v-1) \\\hline
                            \end{array}}\\
                            &=\pm\det\,{\small\begin{array}{|c|c|}
                           \hline W_{\overline{t},k-1}(v-1) & \begin{array}{c} O
                           \\\hline R_{tk}(v-1)\end{array}\\\hline
                            O & W_{t\overline{-}1,k}(v-1) \\\hline
                            \end{array}}\,.\end{align*}
By induction, being square matrices, $W_{\overline{t},k-1}(v-1)$ and
$W_{t\overline{-}1,k}(v-1)$ are both unimodular and so is
$W_{\overline{t}k}(v)$. This completes the proof.}
\end{proof}

 The following corollary is an immediate consequence of Theorem
 \ref{snf}.
\begin{cor}\label{prank}  For any prime $p$, the matrix $W_{\overline{t}k}$ has full $p$-rank.
\end{cor}
\begin{cor}\label{row} Over any field of characteristic not in the set of primes
 $\{p : p|{k-i\choose t-i}, i=0,\ldots,t \}$, the row space of
  $W_{\overline{t}k}$ is the same as that of $W_{tk}$.
\end{cor}
\begin{proof} { Over a field of characteristic not in the set of primes
 $\{p : p|{k-i\choose t-i}, i=0,\ldots,t \}$, the matrix $D_{\overline{t}k}$ is
 non-singular. So the result follows from (\ref{wd}).}\end{proof}

The following corollary was first proved by Wilson \cite{wilson}.
Another proof is also given by Bier \cite{Bier}.

\begin{cor} \label{wil} If  $t\leq k\leq v-t$, then $W_{tk}$ has a diagonal form
which is the ${v\choose t}\times{v\choose k}$ diagonal matrix with
diagonal entries
$${k-i \choose t-i} \hbox{ with multiplicity }{v\choose
i}-{v\choose i-1}, \quad i=0, 1, \ldots, t.$$
\end{cor}
\begin{proof}{ By the proof of Theorem \ref{snf},
$W_{\overline{t}k}=(A\mid B)$, where $A$ is the unimodular
submatrix. We have
$$D_{\overline{t}k}W_{\overline{t}k}=\begin{array}{|c|c|}\hline
D_{\overline{t}k}& O \\\hline\end{array}\,\begin{array}{|c|c|}\hline A& B\\
\hline O & I \\\hline \end{array}\,,$$ where $I$ is the identity
matrix of order ${v\choose k}-{v\choose t}$. Therefore by
(\ref{wd}),
$$W_{\overline{t}t}W_{tk}\, \begin{array}{|c|c|}\hline A^{-1}& -A^{-1}B\\
\hline O & I \\\hline \end{array}=\begin{array}{|c|c|}\hline
D_{\overline{t}k}& O \\\hline\end{array}\,.$$ Since the matrices
$W_{\overline{t}t}$ and $A^{-1}$ are both unimodular, the proof is
completed.}\end{proof}

\begin{remark} We may reformulate the system of (\ref{w1}) in terms of the matrix $W_{\overline{t}k}$.
 By (\ref{wd}), $\bf{x}$ is a solution of (\ref{w1}) if and only if
 $W_{\overline{t}k}{\bf x}=\lambda{\bf h}$, in which $${\bf
h}=D^{-1}_{\overline{t}k}W_{\overline{t}t}{\bf
1}=\left(\tfrac{{v-i\choose t-i}}{{k-i\choose t-i}}^{{v\choose
i}-{v\choose i-1}}, i=0, \ldots,t\right)^{\top},$$ where the
exponents indicate the multiplicity.
\end{remark}

The following corollary was first proved by Wilson in \cite{wil2}.
An alternative proof is given in \cite{wilson}. When applied to a
constant vector it results in necessary and sufficient conditions on
the existence of signed $t$-designs. This case was also proved in
\cite{gju}.

\begin{cor} Let $t\le k\le v-t$. The system
\begin{equation}\label{w2}W_{tk}\bf{x}=\bf{b}\end{equation}
has an integral solution if and only if \begin{equation}\label{b}
\frac{1}{{k-i\choose t-i}}R_{it}\bf{b}\end{equation} is integral for
$i=0,\ldots,t$.
\end{cor}
\begin{proof}{ By (\ref{wd}), ${\bf x}$ is a solution of (\ref{w2}) if and only if
$W_{\overline{t}k}{\bf x}={\bf b'},$
 in which \begin{equation}\label{b'} {\bf
b'}=D^{-1}_{\overline{t}k}W_{\overline{t}t}{\bf b}={\small
\begin{array}{|c|} \hline \frac{1}{{k\choose
t}}R_{0t} \\ \hline \frac{1}{{k-1\choose t-1}}R_{1t}\\
\hline \vdots
\\ \hline \frac{1}{{k-t\choose
t-t}}R_{tt}\\\hline\end{array}}\,\bf{b}\end{equation}
 If (\ref{w2}) has an integral solution, then the right-hand side of
 (\ref{b'}) has to be integral, this proves the necessity.
 Now assume ${\bf b}$ satisfies (\ref{b}). So ${\bf b'}$ is integral.
 By the proof of Theorem \ref{snf}, $W_{\overline{t}k}=(A\mid B)$,
 where $A$ is a unimodular matrix. So its inverse,
$A^{-1}$ is  integral.  Let ${\bf
x}=${\small$\begin{pmatrix}A^{-1}{\bf b'}\\{\bf 0}\end{pmatrix}$}.
Then ${\bf x}$ is integral and $W_{\overline{t}k}{\bf x}={\bf b'}.$
Therefore (\ref{w2}) has an integral solution. }\end{proof}

\begin{remark} Let $M_{tk}(v)$ denote a matrix which is obtained by stacking the
matrices $W_{0k}(v), \ldots, W_{tk}(v)$ one on top of the other. The
rows of $W_{\overline{t}k}$ are in fact a basis for the row space of
$M_{tk}(v)$. The existence of a basis for the row space of
$M_{tk}(v)$ with the property that it contains exactly ${v \choose
i}-{v \choose {i-1}}$ rows from $W_{ik}(v)$, for $i=0,\ldots,t$, was
first observed by Wilson \cite{wilson}. A concrete basis for
$M_{tk}(v)$ which is coincident with $W_{\overline{t}k}$ was
demonstrated by Bier \cite{Bier}.
\end{remark}

\begin{remark} We may represent the result of
Corollary \ref{row} in a more general setting. More precisely if we
replace each row index $T$ of $W_{tk}$  with rank$(T)\le t-1$ by
$T^-$ and keep those indices of rank $t$ and form a new inclusion
matrix, then the resulting matrix is again row-wise
 equivalent to $W_{tk}$. To see this, we fix $t$ and $k$ and simply call
 this new
 matrix $U$. Then  in an appropriate ordering,  $U={\tiny \begin{pmatrix}
 W_{t-1,k}\\
 R_{tk} \\
 \end{pmatrix}}.$ Since $W_{t-1,k}$ and $W_{t\overline{-}1,k}$
 have the same row space, the matrix $U$ has the same row space as
 ${\tiny \begin{pmatrix}
  W_{t\overline{-}1,k}\\
  R_{tk} \\
 \end{pmatrix}}=W_{\overline{t}k}.$
 So by Corollary \ref{row}, we are done . In the same sense
if we replace each row index $T$ of $W_{tk}$ with rank$(T)\le t-m$ by
$T^{-m}$ and each row index $T$ of rank $> t-m$ by $\overline{T}$ and
form a new inclusion matrix, then the resulting matrix is again
row-wise equivalent to $W_{tk}$.
\end{remark}

\begin{remark} It is known that in any decomposition of $2^{[v]}$
into symmetric skipless chains the number of chains beginning at
level $i$ is exactly ${v \choose i}-{v \choose {i-1}}$. One may
guess that the result of Corollary \ref{prank} remains true for any
such decomposition (instead of decomposition of rank chains). But
here is a counterexample. Consider the following decomposition of
the poset $2^{[4]}$ into symmetric skipless chains. {\footnotesize
$$ \emptyset \to 4 \to 14 \to 124 \to 1234,$$
$$ 1 \to 13 \to 134,\ \ \
 2 \to 23 \to 123, \ \ \
 3 \to 24 \to 234,$$
$$12, \ \ \ 34.$$}
If we replace the row indices of $W_{2,2}(4)$ by the first elements
of the chains containing each, we have the following inclusion
matrix which is singular.
{\footnotesize $$\begin{array}{rc} & \begin{array}{cccccc}12& 13& 14& 23& 24& 34\end{array}\\
\begin{array}{c}\emptyset\\ 1 \\2 \\3 \\ 12\\ 34 \end{array}\hspace{-.5 cm}&
\begin{array}{|cccccc|} \hline 1~ &   1 ~ &  1~&    1 ~&   1 ~&   1\\
                       1~ &   1~  &  1~ &   0~ &   0~&    0\\
                       1 ~&   0~&    0  ~ & 1~  &  1 ~ &  0\\
                        0  ~ & 1  ~&  0  ~ & 1~  &  0~   & 1\\
                       1 ~ &  0  ~ & 0 ~ &  0~ &   0~&    0\\
                       0  ~ & 0  ~ & 0  ~ & 0 ~  & 0 ~ &  1\\ \hline
\end{array}\end{array}$$}

\end{remark}

\vspace{5mm} \noindent{\bf\Large 5. The inclusion matrix
$W_{t\underline{k}}(v)$}\vspace{3mm}

In this section we construct another inclusion matrix which is
somehow equivalent to $W_{tk}(v)$ and have the same Smith form as
$W_{\overline{t}k}(v)$. To construct such a matrix, similar to the
previous one, we need a set of (suitable) chains which partition the
poset of $2^{[v]}$. These chains are obtained as follows. In each
chain which is formed by using the notion of rank in Section 3,
replace any subset by its complement. In this way we have new chains
which partition the poset $2^{[v]}$. Now in $W_{tk}$ replace any
column index $K$ by the largest element of the chain containing $K$,
namely $\underline{K}:=[v]\setminus(\overline{[v]\setminus K})$,
keep the indices of the rows and form a new inclusion matrix. We
denote this matrix by $W_{t\underline{k}}(v)$. If $k>v/2$ (resp.
$k\le v/2$), then for each $k$-subset $K$, $k\le|\underline{K}|\le
v$ (resp. $v-k\le|\underline{K}|\le v$). Therefore the minimum size
of the indices of the columns of $W_{t\underline{k}}(v)$ is
$k^*:=\max\{k, v-k\}$.
 Let $Q_{tj}$, $k^*\le j\le v$, be the submatrix of
$W_{t\underline{k}}(v)$, consists of the columns of size $j$. Then
we have the following block decomposition
\begin{equation}\label{wq1}W_{t\underline{k}}=\begin{array}{|c|c|c|c|}
                           \hline  Q_{tk^*} & Q_{t,k^*+1} & \cdots& Q_{tv} \\\hline
                                                      \end{array}\,. \end{equation}
By the same reason as (\ref{rw}), we have
$$W_{it}Q_{tj}= {j-i\choose t-i}Q_{ij} ,\quad\quad i\leq
t\leq k \leq v.$$ Therefore it is seen that
$$W_{it}W_{t\underline{k}}={\small\begin{array}{|c|c|c|c|}
                           \hline {{k^*-i\choose t-i}}Q_{ik^*} & {k^*+1-i\choose t-i}Q_{i,k^*+1}
                           & \cdots& {v-i\choose t-i}Q_{iv} \\\hline
                                                      \end{array}}\,.$$
Note that the number of columns of $Q_{tj}$ is equal to the number
of rank chains whose smallest elements are of size $v-j$. This
number is ${v\choose v-j}-{v\choose v-j-1}={v\choose j}-{v\choose
j+1}$. So it follows that
\begin{equation}\label{wq2}W_{tk}W_{k\underline{k}}=W_{t\underline{k}}D_{t\underline{k}},\end{equation}
where $D_{t\underline{k}}$ is the ${v\choose k}\times{v\choose k}$
diagonal matrix with diagonal entries ${j-t \choose k-t} \hbox{ with
multiplicity }{v\choose j}-{v\choose j+1}, j=k^*, \ldots, v$. Note
that {\small$$\left\{\underline{T} : T\in{[v]\choose
t}\right\}=\left\{[v]\setminus\overline{T} : T\in{[v]\choose
t}\right\}.$$} Let $T_1, T_2\in{[v]\choose t}$, then from the simple
relation
$$T_1\subseteq \underline{T_2}\,\Leftrightarrow\,
[v]\setminus\underline{T_2} \subseteq [v]\setminus T_1,$$ it follows
that
$$W_{t\underline{t}}=\left\{
                       \begin{array}{ll}
                          W_{v\overline{-}t,v-t}^\top, & t> v/2 \hbox{;} \\
W_{\overline{t},v-t}^\top, & t\le v/2.
 \end{array}
                     \right.$$ Being  square matrices,
                      $W_{\overline{t},v-t}$ and $W_{v\overline{-}t,v-t}$ are
unimodular, by Theorem \ref{snf}, and so is $W_{t\underline{t}}$. If
$t\le k \le v-t$, then $k^*\leq t^*$. So it follows from (\ref{wq1})
that $W_{t\underline{t}}$ is a submatrix of $W_{t\underline{k}}$.
Thus $W_{t\underline{k}}$ has a unimodular submatrix of order
${v\choose t}$.

We close our paper with the following summary of the results of this
section.
\begin{thm}\label{last} Let $t\le k \le v-t$. Then the following hold.\\
(i) $W_{t\underline{k}}$ has $(I\mid O)$ as  Smith form, and
consequently has full $p$-rank for any prime $p$.\\
(ii) $W_{t\underline{k}}{\bf x}=\lambda {\bf 1}$ if and only if
$W_{k\underline{k}}D_{t\underline{k}}^{-1}{\bf x}$ is a solution of
(\ref{w1}).
\end{thm}

\noindent{\bf Acknowledgements. } The research of the last author
was in part supported by a grant from University of Tehran.

\end{document}